\documentclass[11pt]{amsart}    
\usepackage{amscd}

\def\cal#1{\mathcal{#1}}

\def\NZQ{\Bbb}               

\def\ZZ{{\NZQ Z}}

\def\frk{\frak}               

\def\mm{{\frk m}}

\def\nn{{\frk n}}

\def\opn#1#2{\def#1{\operatorname{#2}}} 
\opn\chara{char}
\opn\length{\ell}
\opn\pd{pd}
\opn\rk{rk}
\opn\projdim{proj\,dim}
\opn\rank{rank}
\opn\depth{depth}
\opn\grade{grade}
\opn\height{height}
\opn\embdim{emb\,dim}
\opn\codim{codim}

\opn\Tr{Tr}
\opn\bigrank{big\,rank}
\opn\superheight{superheight}\opn\lcm{lcm}
\opn\trdeg{tr\,deg}%
\opn\reg{reg}
\opn\lreg{lreg}
\opn\div{div}
\opn\Div{Div}
\opn\cl{cl}
\opn\Cl{Cl}
\opn\Spec{Spec}
\opn\Proj{Proj}
\opn\Supp{Supp}
\opn\supp{supp}
\opn\Sing{Sing}
\opn\Ass{Ass}
\opn\Ann{Ann}
\opn\Rad{Rad}
\opn\Soc{Soc}
\opn\Ker{Ker}
\opn\Coker{Coker}
\opn\Im{Im}
\opn\Hom{Hom}
\opn\Tor{Tor}
\opn\Ext{Ext}
\opn\End{End}
\opn\Aut{Aut}
\opn\id{id}

\opn\nat{nat}
\opn\pff{pf}
\opn\Pf{Pf}
\opn\GL{GL}
\opn\SL{SL}
\opn\mod{mod}
\opn\ord{ord}
\opn\transpose{\hbox{}^t}
\opn\aff{aff}
\opn\con{conv}
\opn\relint{relint}
\opn\st{st}
\opn\lk{lk}
\opn\cn{cn}
\opn\core{core}
\opn\vol{vol}
\opn\link{link}
\opn\star{star}
\opn\gr{gr}

\def\poly#1#2#3{#1[#2_1,\dots,#2_{#3}]}
\def\pot#1#2{#1[\kern-0.28ex[#2]\kern-0.28ex]}

\opn\dirlim{\underrightarrow{\lim}}
\opn\inivlim{\underleftarrow{\lim}}
\let\dirsum=\oplus
\let\tensor=\otimes
\let\iso=\cong

\let\Sect=\bigcap
\let\Dirsum=\bigoplus

\let\mcone= * 
\let\to=\rightarrow

\let\To=\longrightarrow
\def\Implies{\ifmmode\Longrightarrow \else
     \unskip${}\Longrightarrow{}$\ignorespaces\fi}
\def\implies{\ifmmode\Rightarrow \else
     \unskip${}\Rightarrow{}$\ignorespaces\fi}
\def\iff{\ifmmode\Longleftrightarrow \else
     \unskip${}\Longleftrightarrow{}$\ignorespaces\fi}

\let\:=\colon
\newtheorem{Theorem}{Theorem}
\newtheorem{Lemma}[Theorem]{Lemma}
\newtheorem{Corollary}[Theorem]{Corollary}
\newtheorem{Proposition}[Theorem]{Proposition}

\let\epsilon\varepsilon
\let\kappa=\varkappa
\opn\inii{in}
\opn\inim{inm}
\opn\set{set}
\def\pnt{{\raise0.5mm\hbox{\large\bf.}}}

\begin{document}

\title{Ideals whose associated graded rings are\\ isomorphic to
the base rings}
\author{Yukihide Takayama}
\address{Yukihide Takayama, Department of Mathematical
Sciences, Ritsumeikan University, 
1-1-1 Nojihigashi, Kusatsu, Shiga 525-8577, Japan}
\email{takayama@se.ritsumei.ac.jp}

\def\Coh#1#2{H_{\mm}^{#1}(#2)}
\def\gCoh#1#2#3{H_{#1}^{#2}\left(#3\right)}

\newcommand{\AppTh}{Theorem~\ref{approxtheorem} }
\def\da{\downarrow}
\newcommand{\ua}{\uparrow}
\newcommand{\namedto}[1]{\buildrel\mbox{$#1$}\over\rightarrow}
\newcommand{\bdel}{\bar\partial}
\newcommand{\proj}{{\rm proj.}}

\maketitle

\newenvironment{myremark}[1]{{\bf Note:\ } \dotfill\\ \it{#1}}{\\ \dotfill
{\bf Note end.}}
\newcommand{\transdeg}[2]{{\rm trans. deg}_{#1}(#2)}
\newcommand{\mSpec}[1]{{\rm m\hbox{-}Spec}(#1)}
\def\tbf{{\Large To Be Filled!!}}

\pagestyle{plain}

\def\nzqK{{\NZQ K}}

\begin{abstract}
Let $k$ be a field.
We determine the ideals $I$ in a finitely generated graded $k$-algebra $A$,
whose associated graded rings $\Dirsum_{n\geq 0}I^n/I^{n+1}$ are
isomorphic to $A$. Also we compute the graded local cohomologies of 
the Rees rings $A[I t]$ and give the condition for $A[I t]$ to be 
generalized Cohen-Macaulay under the condition that $A$ is 
generalized Cohen-Macaulay.
\end{abstract}

\section*{Introduction}
Let $k$ be a field and $S = \poly{k}{X}{n}$ be a polynomial ring over
$k$. Let $J \subset S$ be a homogeneous ideal generated in
degree $\geq 2$ and consider the finitely generated
standard graded $k$-algebra $A = S/J$. 
Let $I = (\bar{f}_1, \ldots, \bar{f}_m) \subset A$ be a homogeneous
ideal and, for $f\in S$, we denote by $\bar{f}$ the image of
$f$ in $A$. In particular, we denote $\bar{X}_i$ by $x_i$.
We consider the associated graded ring $G = \Dirsum_{\ell=0}^\infty
I^\ell/I^{\ell + 1}$ and the Rees ring 
$R = \Dirsum_{\ell=0}^\infty I^\ell t^\ell = A[It]$ with 
regard to $I \subset A$. We set $R_+ = \Dirsum_{\ell=1}^\infty 
I^\ell t^\ell$.

It is well known that $A \iso G$ if $I = \mm := (x_1,\ldots, x_n)$.  In
this paper, we determine precisely the ideals $I \subset A$ satisfying
$A\iso G$. Such ideals turn out to be generated by certain linear
forms (Theorem~\ref{main}). Next we compute graded local cohomologies
of Rees algebra with regard to such papameter ideals, which is an 
extension of the result in \cite{HPT2} for $I=\mm$ 
(Proposition~\ref{localcohomology}). As an application, we 
give the condition for generalized Cohen-Macaulayness of the Rees
algebra when the base ring is also generalized Cohen-Macaulay 
(Theorem~\ref{genCMcondition}).

\section{The Main Theorem }
Consider the following natural surjective homomorphism
from the polynomial ring 
$k[\overline{X},\overline{Y}] = k[X_1,\ldots, X_n, Y_1,\ldots, Y_m]
( = S[Y_1,\ldots, Y_m])$
in the variables $X_1,\ldots, X_n, Y_1,\ldots, Y_m$:
\begin{equation}
\label{theMap}
\begin{array}{cccc}
\psi\; : & k[X_1,\ldots, X_n, Y_1,\ldots, Y_m]  & \To & G  \\
         &     X_i             &\longmapsto & [x_i] \in A/I\\ 
         &     Y_i             &\longmapsto & [\bar{f}_i] \in I/I^2 
\end{array}
\end{equation}
where $[a]$ denotes the equivalent class of $a\in A$ in $A/I$ or $I/I^2$. 
For an element $g\in k[X_1,\ldots, X_n, Y_1,\ldots, Y_m]$,
we will denote by $\deg_{\overline{Y}}(g)$ the degree of $g$ 
with regard 
to the variables $Y_1,\ldots, Y_m$.
$\Ker\psi$ is computed as follows.
\begin{Lemma}
\label{theKernel}
$\Ker\psi = J + (f_1,\ldots, f_m) + L$,
where $L$ is an ideal generated by the following set:
\begin{equation*}
\left\{
        g \in k[\overline{X}, \overline{Y}]\; : \; 
	\begin{array}{l}
	 g\notin J + (f_1,\ldots, f_m),\\
        g \mbox{ is homogeneous in $\overline{X}$ and 
	 $\overline{Y}$ respectively},\\
	g(X_1,\ldots, X_n, f_1,\ldots, f_m) 
	\in (f_1,\ldots, f_m)^{d+1} + J\\
        \mbox{with }d = \deg_{\bar{Y}}(g) (>0)
	\end{array}
\right\}
\end{equation*}
\end{Lemma}
\begin{proof}
Consider the natural surjection:
\begin{equation}
\begin{array}{cccc}
\varphi\; : & A[Y_1,\ldots, Y_m]  & \To & G  \\
            &     a \; (\in A)    &\longmapsto & [a] \in A/I \\
            &     Y_i             &\longmapsto & [\bar{f}_i] \in I/I^2 
\end{array}
\end{equation}
where we know that $IA[Y_1,\ldots, Y_m]\subset \Ker\varphi$.
From this we obtain the following commutative diagram:
\begin{equation*}
\begin{CD}
     k[X_1,\ldots, X_n,Y_1,\ldots, Y_m]  @>{\psi}>> G  \\
     @V{\mod J S[\overline{Y}]}VV       @AA{\{Y_i \mapsto [\bar{f}_i]\}_{i=1}^m}A \\
    A[Y_1,\ldots, Y_m] @>{\mod I A[\overline{Y}]}>> (A/I)[Y_1,\ldots, Y_m]           \\
     @.                                         ||  \\
     @. \displaystyle{\frac{k[X_1,\ldots, X_n, Y_1,\ldots, Y_m]}
                      {(J + (f_1,\ldots, f_m))}}
\end{CD}
\end{equation*}
and we know that
 $\psi(g(\overline{X},\overline{Y}))=g(\overline{X}, f_1,\ldots, f_m)
+ (f_1,\ldots, f_m)^{d+1} + J$ if $d=\deg_{\overline{Y}}(g)$.
Then we know that $\Ker\psi$ is generated by the set
\begin{eqnarray*}
{\cal L}' = \{g \in k[\overline{X},\overline{Y}]\; 
  : \;   g(\overline{X}, f_1,\ldots, f_m) 
\in (f_1,\ldots, f_m)^{d+1} + J 
\mbox{ with } d = \deg_{\bar{Y}}g\}.
\end{eqnarray*}
We also know that $J + (f_1,\ldots, f_m) \subset \Ker\psi$.
Since the generators of $J + (f_1,\ldots, f_m)$ 
do not contain $Y_1, \ldots, Y_m$, these generators
are contained in ${\cal L}'$.
\end{proof}

Now we show our main theorem.

\begin{Theorem}
\label{main}
$A \iso G$ if and only if $I$ is generated by  
$x_{i_1},\ldots, x_{i_m}$ such that 
we can choose a set $\{h_1,\ldots, h_p\}$ of 
homogeneous generators of $J$ such that, for
each  $j$, $h_j \in k[X_{i_1},\ldots, X_{i_m}]$
or $h_j \in k[X_{i_{m+1}},\ldots, X_{i_n}]$, 
where $\{i_1,\ldots, i_{m},i_{m+1},\ldots, i_n\} = \{1,\ldots, n\}$.
\end{Theorem}
\begin{proof} We first prove the only-if part. 
We have
\begin{equation*}
G \iso k[X_1,\ldots, X_n, Y_1, \ldots, Y_m]/(J + (f_1,\ldots, f_m) + L)
\end{equation*}
where $L$ is the ideal as given in Lemma~\ref{theKernel}.
Now $J$ and $(f_1,\ldots, f_m)$ give relations on
the variables $X_1,\ldots, X_n$ and $f_i\notin J$ for all $i=1,\ldots, m$.
Thus in order for $G$ to be isomorphic to $A = k[X_1,\ldots, X_n]/J$,
we must have
\begin{enumerate}
\item [$(i)$] $f_1,\ldots, f_m$ are  linear forms in $X_1,\ldots, X_n$, 
say $f_1 = X_1, \ldots, f_m = X_m$ for simplicity, and 
\item [$(ii)$] $L$ precisely gives the relations on $Y_1,\ldots Y_m$, that,
when $Y_i$ are replaced by $X_i$ ($i=1,\ldots, m$), give
the relations on $X_1, \ldots, X_m$ given by $J$. 
\end{enumerate}
In fact, if we have a non-linear relation $f_j$, then $L$ must contain
relations such as $\{X_i-Y_i\}_i\cup \{Y_i\}_i$ to remove the extra
relations on $X_1,\ldots, X_n$ caused by $f_j$. But such relations cannot
be in $L$. Thus we know $(i)$, and from this we also know $(ii)$.

Now by $(i)$, we can assume 
without loss of generality that $I = (x_1,\ldots, x_m)$ for some $m\leq n$,
i.e.,  $f_i = X_i$, $i=1,\ldots, m$. Then
we have 
$G \iso k[X_1,\ldots, X_n, Y_1,\ldots, Y_m]/(J + (X_1,\ldots, X_m) + L)$ 
where $L$ is generated by the set 
\begin{equation*}
{\cal L} :=
\left\{
   g \in k[\overline{X},\overline{Y}] 
  \;:\;
  \begin{array}{l}
   g(\overline{X},\overline{Y}) 
   \mbox{ is homogeneous in $\overline{X}$ and $\overline{Y}$
   respectively } \\
   g\notin J + (X_1,\ldots, X_m),\;  \deg_{\overline{Y}}(g) = d >0, \\
   g(\overline{X}, X_1,\ldots, X_m)\in (X_1,\ldots, X_m)^{d+1} + J
  \end{array}
\right\}
\end{equation*}
For an element $g\in {\cal L}$,
if $g(\overline{X}, X_1,\ldots, X_m) \in (X_1,\ldots, X_m)^{d+1}$
with $d = \deg_{\overline{Y}}(g)$ 
then $g\in (X_1,\ldots, X_m)$, which is a contradition 
since $g\notin (X_1,\ldots, X_m) + J$. 
This we can assume that  $g(\overline{X}, X_1,\ldots, X_m)\in J$.
Thus 
\begin{equation*}
{\cal L} = 
\left\{
   g \in k[\overline{X},\overline{Y}] 
  \;:\;
  \begin{array}{l}
   g \mbox{ is homogeneous in $\overline{X}$ and $\overline{Y}$ respectively},\\
   g\notin J + (X_1,\ldots, X_m),
   g(\overline{X}, X_1,\ldots, X_m)\in  J
  \end{array}
\right\}
\end{equation*} and 
\begin{eqnarray*}
G & \iso  & 
\frac{k[\overline{X},\overline{Y}]/(X_1,\ldots, X_m)}
     {(J + L + (X_1,\ldots, X_m)) /  (X_1,\ldots, X_m)} \\
  & \iso &
   \frac{k[Y_1,\ldots, Y_m, X_{m+1},\ldots, X_n]}
     {(J + L + (X_1,\ldots, X_m)) /  (X_1,\ldots, X_m)}.
\end{eqnarray*}
Now by $(ii)$, we must have 
(a) we can choose a set of generators $h_1,\ldots, h_p$ of $J$
as follows: $h_i \in (X_1,\ldots, X_m)$ or 
$h_i$ contains no monomial from $(X_1,\ldots, X_m)$
for all $i$, 
and 
(b) $L = \sigma(J\cap (X_1,\ldots, X_m))k[Y_1,\ldots, Y_m, X_{m+1},\ldots, X_n]$ 
with $\sigma : k[X_1,\ldots, X_n]\to k[Y_1,\ldots, Y_m, X_{m+1},\ldots, X_n]$
such that $\sigma(X_i) = Y_i$ for $1\leq i\leq m$ and $\sigma(X_i) = X_i$ for
$m+1\leq i \leq n$.
But by the condition (a), ${\cal L}$ can be chosen as follows:
\begin{equation*}
{\cal L} = 
\left\{
   g \in k[\overline{Y}] 
  \;:\;
  \begin{array}{l}
   g(\overline{Y}) \mbox{ is homogeneous},\; 
   g(X_1,\ldots, X_m)\in  J
  \end{array}
\right\}
\end{equation*}
and the condition (b) is satisfied.
Consequently, we must have the following:
we can choose a set of homogeneous generators $h_1,\ldots, h_p$ of $J$
such that as follows: $h_i \in k[X_1,\ldots, X_m]$ or 
$h_i \in k[X_{m+1},\ldots, X_n]$
for each $i$.

The proof of the if-part is carried out by 
tracing the above discussion conversely.
\end{proof}

Now we have the following well-known result as a corollary.

\begin{Corollary}
If $I = \mm$, then $A\iso G$.
\end{Corollary}

\section{Local cohomologies of Rees algebras}

We compute here the local cohomology of our Rees algebra $R$
in the case of $A\iso G$. 
For a graded ring $R = \Dirsum_{n\geq 0}R_n$ such that $R_0$  is local 
with the maximal ideal $\mm$ and a graded $R$-module $M$, we denote
by $\gCoh{{\cal M}}{q}{M}$ the $q$th graded local cohomology,
where ${\cal M} = \mm \dirsum \Dirsum_{n\geq 1} R_n$.

Now we cite two results.

\begin{Theorem}[Herzog-Popescu-Trung~\cite{HPT2}]
\label{HPT2:Theorem3.2}
Let $A = k[X_1,\ldots, X_n]/J$ be a residue class ring with 
regard to a homogeneous ideal $J$ and let $\mm = (x_1,
\ldots, x_n)$. Then we have 
\[
\gCoh{\mm}{i}{A[\mm t]}_a
= \left\{
   \begin{array}{ll}
    \gCoh{\mm}{i}{A}_a^{\dirsum^{a+1}}  & a\geq 0 \\
    0                                                  & a = -1  \\
    \gCoh{\mm}{i-1}{A}_a^{\dirsum^{-a-1}}  & a\leq -2 \\
   \end{array}
  \right.
\]
\end{Theorem}

\begin{Theorem}[Goto-Watanabe~\cite{GW}]
\label{GW:Theorem2.2.5}
Let $R$, $S$ be graded rings defined over $k$ and 
$\mm = R_+$, $\nn = S_+$ be their H-maximal ideals.
We put $T = R\tensor_k S$ and ${\cal M} = T_+$.
If $M$ (resp. $N$) is a graded $R$- (resp. $S$-) module,
we have
\begin{equation*}
    \gCoh{{\cal M}}{q}{M\tensor_k N}
= \Dirsum_{i+j=q} \gCoh{\mm}{i}{M} \tensor_k \gCoh{\nn}{j}{N}
\end{equation*}
\end{Theorem}

Now, according to Theorem~\ref{main}, we can assume that 
$I = (x_{i_1},\ldots, x_{i_m}) \subset A$ and 
$A = S/J$ with $J=(h_1,\ldots, h_r, h_{r+1}, \ldots, h_p)$
such that $h_1,\ldots, h_r \in k[X_{i_1},\ldots, X_{i_m}]$
and $h_{r+1},\ldots, h_p\in k[X_{i_{m+1}},\ldots, X_{i_n}]$.
We further set $A_1 = k[X_{i_1},\ldots, X_{i_m}]/(h_1,\ldots, h_r)$,
$\mm_1 = (x_{i_1},\ldots, x_{i_m}) \subset A_1$
and $A_2 = k[X_{m+1},\ldots, X_n]/(h_{r+1},\ldots, h_p)$,
$\mm_2 = (x_{i_{m+1}},\ldots, x_{i_n}) \subset A_2$.
Then we have $A = A_1\tensor_k A_2$.
\begin{Proposition} 
\label{localcohomology}
Let $A$ and $J$ be as above.
Let $R = A[I t]$ be the Rees algebra with regard to
$I = (x_{i_1},\ldots, x_{i_m}) \subset A$ where
$X_{i_1}, \ldots, X_{i_m}$ are as in Theorem~\ref{main}.
Then we have 
\begin{eqnarray*}
\label{cohomology}
\gCoh{{\cal M}}{\ell}{R}_a 
  & = & 
        \Dirsum_{\tiny\begin{array}{c}
	 \alpha + \beta = a\\
         \alpha\geq 0
	  \end{array}}
        \Dirsum_{p+q=\ell}
                       \gCoh{\mm_1}{p}{A_1}_\alpha^{\dirsum^{\alpha+1}}
                       \tensor_k \gCoh{\mm_2}{q}{A_2}_\beta \\
  &   &\dirsum
        \Dirsum_{\tiny\begin{array}{c}
	 \alpha +\beta= a\\
         \alpha\leq -2
	  \end{array}}
        \Dirsum_{p+q=\ell}
	\gCoh{\mm_1}{p-1}{A_1}_\alpha^{\dirsum^{-\alpha-1}}
                       \tensor_k \gCoh{\mm_2}{q}{A_2}_\beta \\
\end{eqnarray*}
where ${\cal M} = \mm\tensor R_+$.
\end{Proposition}

\begin{proof}
We have $A[ I t] \iso A_1[\mm_1 t]\tensor_k A_2$.
Thus by Theorem~\ref{GW:Theorem2.2.5} and  Theorem~\ref{HPT2:Theorem3.2}
we compute
\begin{eqnarray*}
\gCoh{{\cal M}}{\ell}{A[It]}_a
  & = & \Dirsum_{i+j=\ell} (\gCoh{\mm_1}{i}{A_1[\mm_1 t]}
                       \tensor_k \gCoh{\mm_2}{j}{A_2})_a \\
  & = & \Dirsum_{i+j=\ell}
        \Dirsum_{\alpha + \beta = a}
  	  \gCoh{\mm_1}{i}{A_1[\mm_1 t]}_\alpha
                       \tensor_k \gCoh{\mm_2}{j}{A_2}_\beta \\
  & = & 
        \Dirsum_{\tiny\begin{array}{c}
	 \alpha + \beta = a\\
         \alpha\geq 0
	  \end{array}}
        \Dirsum_{i+j=\ell}
                       \gCoh{\mm_1}{i}{A_1}_\alpha^{\dirsum^{\alpha+1}}
                       \tensor_k \gCoh{\mm_2}{j}{A_2}_\beta \\
  &   &\dirsum
        \Dirsum_{\tiny\begin{array}{c}
	 \alpha + \beta = a\\
         \alpha\leq -2
	  \end{array}}
        \Dirsum_{i+j=\ell}
	\gCoh{\mm_1}{i-1}{A_1}_\alpha^{\dirsum^{-\alpha-1}}
                       \tensor_k \gCoh{\mm_2}{j}{A_2}_\beta \\
\end{eqnarray*}
as required.
\end{proof}

We will
denote the $a$-invariant of a graded $k$-algebra $R$
by $a(R) = \max\{j : \gCoh{\mm}{\dim R}{R}_j\ne 0\}$.

Notice that, since $A\iso G$, Cohen-Macaulayness of the Rees algebra
is trivial by Theorem~5.1.22~\cite{TI}. Namely, $R$ is Cohen-Macaulay
if and only if $A$ is Cohen-Macaulay and $a(A) < 0$. 

Recall that a finitely generated graded $k$-algebra $R$ 
is a generalized Cohen-Macaulay ring if $\ell(\gCoh{R_+}{i}{R}) <
\infty$ for all $i< \dim R$, where $\ell(-)$ denotes the length.
Now we consider
generalized Cohen-Macaulayness of the Rees algebra under the condition
that $A$ is generalized Cohen-Macaulay. 
If $I = \mm$, this is immediate 
from Theorem~\ref{HPT2:Theorem3.2}. Namely, if $A$ is generalized 
Cohen-Macaulay, so is $R$. We now consider the general case.
Before that we prepare a lemma.

\begin{Lemma}
\label{genCM:tensor}
Let $R_1$ and $R_2$ be finitely generated  graded $k$-algebra. Then 
if $R_1\tensor_k R_2$ is generalized Cohen-Macaulay ring
then so are $R_i$ $(i=1,2)$.
\end{Lemma}
\begin{proof}
Let $B := R_1\tensor_k R_2$, $\mm = R_+$, $\mm_i := (R_i)_+$ ($i=1,2$)
and $d := \dim B$. Then we have by Theorem~\ref{GW:Theorem2.2.5} 
\begin{equation*}
    0 = \gCoh{\mm}{i}{B}_P
      = \Dirsum_{n_1+n_2 = i} 
	(\gCoh{\mm_1}{n_1}{R_1}\tensor_k\gCoh{\mm_2}{n_2}{R_2})_P
\end{equation*}
for arbitrary $i\ne d$ and 
$P\in\Spec{B}\backslash\{\mm\}$. Since $P = P_1\tensor_k R_2
+ R_1\tensor_k P_2$ for some primes $P_1$ and $P_2$ such that 
$(P_1, P_2)\ne (\mm_1, \mm_2)$, we have 
$(\gCoh{\mm_1}{n_1}{R_1}\tensor_k\gCoh{\mm_2}{n_2}{R_2})_P
  \iso \gCoh{\mm_1}{n_1}{R_1}_{P_1}\tensor_k\gCoh{\mm_2}{n_2}{R_2}_{P_2}$.
Thus we have 
\begin{equation*}
\gCoh{\mm_1}{n_1}{R_1}_{P_1}\tensor_k\gCoh{\mm_2}{n_2}{R_2}_{P_2} = 0
\end{equation*}
for aribtrary $(P_1, P_2)\ne (\mm_1, \mm_2)$ and $n_1 + n_2 \neq d$.
From this we  know that $\gCoh{\mm_i}{n_i}{R_i}_{P_i} = 0$
for arbitrary $n_i \ne \dim R_i$ and prime ideal $P_i\ne \mm_i$
($i=1,2$).
\end{proof}

Now we give a condition for $R$ such that 
$A\iso G$ to be generalized Cohen-Macaulay.

\begin{Theorem}
\label{genCMcondition}
Let $A = k[X_1,\ldots, X_n]/J$ be a graded generalized Cohen-Macaulay
ring and $I \subset A$ be such that $A\iso G$. Also let 
$P_1,\ldots, P_u$ be the minimal primes of $A$ such that 
$\dim A = \dim A/P_i$.
Then the Rees ring
$R = A[It]$ is generalized Cohen-Macaulay if and only if 
one of the following holds:
\begin{enumerate}
\item $I \subset \Sect_{i=1}^u P_i$
\item $I \not\subset \Sect_{i=1}^u P_i$ and $\dim A_2 =0$

\item $I \not\subset \Sect_{i=1}^u P_i$, 
$\dim A_1 >0$, $\dim A_2 >0$, 
$a(A_1)<0$, $\gCoh{\mm_1}{d_1-1}{A_1}_j=0$ for all $j\leq -2$,
and $\gCoh{\mm_2}{d_2-1}{A_2}=0$.
\end{enumerate}
\end{Theorem}
\begin{proof}
We have 
\begin{equation*}
\dim R = \left\{
	  \begin{array}{ll}
	   \dim A      & \mbox{if } I 
 	                 \subset \Sect_{i=1}^u P_i \\
	   \dim A +1   & \mbox{otherwise}
	  \end{array}
	 \right.
\end{equation*}
where $P_1,\ldots, P_u$ are the minimal prime ideals of $A$
such that $d = \dim A/P_i$ (see \cite{V}).

Now from the short exact sequences 
\begin{equation*}
 0 \To R_+ \To R \To A \To 0
\qquad\hbox{and}\qquad
 0 \To R_+(1)\To R \To G \To 0,
\end{equation*}
we have, by assumption, the exact sequences
\begin{equation*}
0= \gCoh{\mm}{i-1}{A}_p 
\To \gCoh{{\cal M}}{i}{R_+}_P
\To \gCoh{{\cal M}}{i}{R}_P
\To \gCoh{\mm}{i}{A}_p =0
\end{equation*}
and 
\begin{equation*}
 0= \gCoh{{\cal M}}{i-1}{G}_P 
\To \gCoh{{\cal M}}{i}{R_+}_P(1)
\To \gCoh{{\cal M}}{i}{R}_P
\To \gCoh{{\cal M}}{i}{G}_P =0
\end{equation*}
for $i\ne d, d+1$ ($d:=\dim A$) and 
$P\in\Spec{R}\backslash \{{\cal M}\}$ with
$p = P\cap A$. 
Thus $\gCoh{{\cal M}}{i}{R}_P(1) \iso \gCoh{{\cal M}}{i}{R}_P$
for such $P$ and $i$. 
Then we have $\gCoh{{\cal M}}{i}{R}_P = \gCoh{{\cal M}}{i}{R}_P(1)$
for $i\ne d, d+1$ and $P \ne {\cal M}$. Then, since $\gCoh{{\cal M}}{i}{R}$
is Artinian, we have $\gCoh{{\cal M}}{i}{R}_P=0$, i.e., 
$\length(\gCoh{{\cal M}}{i}{R}) < \infty$
for $i=0,\ldots, d-1$. Thus $R$ is generalized Cohen-Macaulay if $\dim R = \dim A$.

We know consider the case of $\dim R = \dim A + 1$.
We set $d_i := \dim A_i$ ($i=1,2$) and we then  have $d = d_1 + d_2$.
By Proposition~\ref{localcohomology}, we have
\begin{eqnarray*}
\lefteqn{\gCoh{{\cal M}}{d}{R}_a }\\
  & = & 
        \Dirsum_{\tiny\begin{array}{c}
	 \alpha + \beta = a\\
         \alpha\geq 0
	  \end{array}}
                       \gCoh{\mm_1}{d_1}{A_1}_\alpha^{\dirsum^{\alpha+1}}
                       \tensor_k \gCoh{\mm_2}{d_2}{A_2}_\beta \\
  &   &\dirsum
        \Dirsum_{\tiny\begin{array}{c}
	 \alpha +\beta= a\\
         \alpha\leq -2
	  \end{array}}
	\left[
	\gCoh{\mm_1}{d_1-1}{A_1}_\alpha^{\dirsum^{-\alpha-1}}
                       \tensor_k \gCoh{\mm_2}{d_2}{A_2}_\beta 
	\dirsum 
        \gCoh{\mm_1}{d_1}{A_1}_\alpha^{\dirsum^{-\alpha-1}}
                       \tensor_k \gCoh{\mm_2}{d_2-1}{A_2}_\beta 
	\right]
\end{eqnarray*}
so that 
\begin{eqnarray*}
\lefteqn{\gCoh{{\cal M}}{d}{R}}\\
 & = & \Dirsum_{a\in\ZZ} \gCoh{{\cal M}}{d}{R}_a \\
 & = & 
        \Dirsum_{\alpha\geq 0}
            \gCoh{\mm_1}{d_1}{A_1}_\alpha^{\dirsum^{\alpha+1}}
            \tensor_k \gCoh{\mm_2}{d_2}{A_2}
      \dirsum
        \Dirsum_{\alpha\leq -2}
   	    \gCoh{\mm_1}{d_1-1}{A_1}_\alpha^{\dirsum^{-\alpha-1}}
            \tensor_k \gCoh{\mm_2}{d_2}{A_2}\\
 &    &
       \dirsum 
        \Dirsum_{\alpha\leq -2}
            \gCoh{\mm_1}{d_1}{A_1}_\alpha^{\dirsum^{-\alpha-1}}
            \tensor_k \gCoh{\mm_2}{d_2-1}{A_2}.
\end{eqnarray*}
Now by Lemma~\ref{genCM:tensor}, 
$\gCoh{\mm_1}{d_1-1}{A_1}$ and $\gCoh{\mm_2}{d_2-1}{A_2}$
are of finite length. Also 
$\gCoh{\mm_1}{d_1}{A_1}$ and  $\gCoh{\mm_2}{d_2}{A_2}$
are non-zero Artinian modules. 
Thus if $\length(\gCoh{\mm_2}{d_2}{A_2})=\infty$, 
we know that we have $\length(\gCoh{{\cal M}}{d}{R}) <\infty$ provided
\begin{enumerate}
\item $\gCoh{\mm_2}{d_2 -1}{A_2}=0$
      or $\length(\gCoh{\mm_1}{d_1}{A_1})<\infty$,
\item $a(A_1)< 0$ and
\item $\gCoh{\mm_1}{d_1 -1}{A_1}_j=0$ for all $j\leq -2$.
\end{enumerate}
On the other hand, if 
$\length(\gCoh{\mm_2}{d_2}{A_2})<\infty$,
we have $\length(\gCoh{{\cal M}}{d}{R}) <\infty$ provided
$\gCoh{\mm_2}{d_1 -1}{A_2}=0$ or
 $\length(\gCoh{\mm_1}{d_1}{A_1})<\infty$.
By Grothendieck's finiteness theorem (see, for example, Theorem~9.5.2
\cite{BS}), we have 
$\length(\gCoh{\mm_i}{d_i}{A_i})<\infty$
if and only if $\dim A_i = 0$ ($i=1,2$). Thus we have the 
desired result.
\end{proof}


\end{document}